\documentclass{article}%
\usepackage{amsmath}
\usepackage{amsfonts}
\usepackage{amssymb}
\usepackage{graphicx}
\usepackage{subfigure}
\usepackage{epsfig}%
\usepackage{color}
\usepackage{multirow}
\usepackage{setspace}
\usepackage{hyperref}
\newtheorem{theorem}{Theorem}

\newtheorem{lemma}[theorem]{Lemma}

\newtheorem{proposition}[theorem]{Proposition}

\newenvironment{proof}[1][Proof]{\noindent\textbf{#1.} }{\ \rule{0.5em}{0.5em}}
\begin{document}
\doublespacing
\title{The Two-Sample Problem Via Relative Belief Ratio}

\author{Luai Al-Labadi\thanks{{\em Address for correspondence:} Luai Al-Labadi, Department  of Mathematics, University of Sharjah, P. O. Box 27272, Sharjah, UAE. E-mail:  lallabadi@sharjah.ac.ae.}}



\date{}
\maketitle

\pagestyle {myheadings} \markboth {} {Two-Sample Problem: RB Ratio}

\begin{abstract}
This paper deals with a new Bayesian approach to the two-sample problem. More specifically, let $x=(x_1,\ldots,x_{n_1})$  and $y=(y_1,\ldots,y_{n_2})$ be two independent samples coming from unknown distributions  $F$ and $G$, respectively.  The goal is to test the null hypothesis $\mathcal{H}_0:~F=G$ against all possible alternatives.   First, a Dirichlet process prior for $F$ and $G$ is considered. Then the change of their Cram\'{e}r-von Mises distance from a priori to a posteriori is compared through the relative belief ratio. Many theoretical properties of the procedure have been developed and several examples have been discussed, in which the proposed approach shows excellent performance.
\par

 \vspace{9pt} \noindent\textsc{Keywords:}  Dirichlet process,  hypothesis testing,  relative belief inferences, two-sample problem.

 \vspace{9pt}


\end{abstract}

\section{Introduction}

For two independent samples, the \emph{two-sample problem} is concerned to determine whether the two samples are generated from the same population. Although it is considered an old problem in statistics, it always attracts the attention of researchers  due to it applications in different fields. For instance, in medical studies, one may want to asses the efficiency of a new drug to two groups of patients.

The two-sample problem can be stated formally as follows. Given two independent samples $x=(x_1,\ldots,x_{n_1}) \overset {i.i.d.} \sim F$ and $y=(y_1,\ldots,y_{n_2}) \overset {i.i.d.} \sim G$, with $F$ and $G$ being unknown continuous cumulative distribution functions  (cdf's), the aim  is to test the null hypothesis  $\mathcal{H}_0:~F=G$ against all other alternatives.

The methodology developed in this paper is Bayesian and it is inspired from the recent work of Al-Labadi and Evans (2018) for model checking.  At first, two Dirichlet processes $DP(a_1, H_1)$ and $DP(a_2,H_2)$ are  considered as  priors for $F$ and $G$, respectively.  Then the concentration of the posterior distribution of the distance between the two processes is compared to the concentration of the prior distribution of the distance between the two processes.
If the posterior is more concentrated about the model than the prior, then this is evidence in favor of $\mathcal{H}_0$ and if the posterior is
less concentrated, then this is evidence against $\mathcal{H}_0$. This comparison is made through a particular measure of evidence known as the \emph{relative belief ratio}, which will indicate  whether there is evidence for or against $\mathcal{H}_0$.  Moreover, a calibration of this evidence is provided concerning whether there is strong or weak evidence for or against the hypothesis. The proposed methodology is simple, general and does not require obtaining a closed form of the relative belief ratio. More details about relative belief ratio are highlighted in Section 2 of this paper.

Developing procedures for hypothesis testing has recently given a considerable attention in the literature of Bayesian nonparametric inference.  A main stream of these procedures has focused on embedding the suggested model as a null hypothesis in a larger family of distributions. Then priors are placed on the null and the alternative and a Bayes factor is computed.  For instance, Florens, Richard, and Rolin (1996) used a Dirichlet process for the prior on the alternative.  Carota and Parmigiani (1996), Verdinelli and Wasserman (1998), Berger and Guglielmi (2001) and McVinish, Rousseau, and Mengersen (2009) considered a mixture of Dirichlet processes, a mixture of Gaussian processes, a mixture of P\'{o}lya trees and a mixture of triangular distributions, respectively, for the prior on the alternative. Another approach for model testing is based on placing a prior on the true distribution generating the data and measuring the distance between the posterior distribution and the proposed one. Swartz (1999) and Al-Labadi and Zarepour (2013, 2014a) considered the Dirichlet process prior and used the Kolmogorov distance to derive a goodness-of-fit test for continuous models. Viele (2000) used the Dirichlet process and the Kullback-Leibler distance to test discrete models. Hsieh (2011) used the P\'{o}lya tree prior and the Kullback-Leibler distance to test continuous distributions. The work described above focuses only on goodness of fit tests and model checking. With regard to the two-sample problem, the literature is very scarce and scattered. Some exceptions include the remarkable work of Holmes, Caron, Griffin, and Stephens (2015) who developed a way to compute the Bayes factor for testing the null hypothesis through the marginal likelihood of the data with P\'olya tree priors centered either subjectively or using an empirical procedure. Under the null hypothesis, they modeled the two samples to come from a single random measure distributed as a P\'olya tree, whereas under the alternative hypothesis the two samples come from two separate P\'olya tree random measures. Ma and Wong (2011) allowed the two distributions to be generated jointly through optional coupling of a P\'olya tree prior. Borgwardt and Ghahramani (2009) discussed two-sample tests based on Dirichlet process mixture models and derived a formula to compute the Bayes factor in this case. An extension of the Bayes factor approach based on P\'olya tree priors to cover censored and multivariate data was proposed by Chen and Hanson (2014). Huang and Ghosh (2014) considered the two-sample hypothesis testing problems under P\'olya tree priors and Lehmann alternatives. Shang and Reilly (2017) introduced a class of tests, which use the connection between the Dirichlet process prior and the Wilcoxon rank sum test. They also extend their idea using the Dirichlet process mixture prior and developed a Bayesian counterpart to the Wilcoxon rank sum statistic and the weighted log rank statistic for right and interval censored data. In a recent work, Al-Labadi and Zarepour (2017) proposed a method based on the Kolmogorov distance and  samples from the Dirichlet process  to assess the equality of two unknown distributions, where the distance between two posterior Dirichlet processes  is compared with a reference  distance. The parameters of the two  Dirichlet processes are chosen so that    any  discrepancy between   the  posterior distance  and the reference distance is only attributed to the difference between the two samples.

In Section 3, the Dirichlet process prior $DP(a,H)$ is briefly reviewed. In Section 4,  the Cram\'{e}r-von Mises distance between two Dirichlet processes is considered and several  of its theoretical  properties are developed. Section 5 addresses setting parameters of the two Dirichlet processes. In Section 6, a computational algorithm of the approach is developed. Section 7 presents several examples where the behaviour of the approach is inspected. Finally, some concluding remarks are made in Section 8. The proofs are placed in the Appendix.

\section{Relative Belief Ratios}

In this section, for the reader's convenience,  some background of  relative belief ratios is provided.  For more details about this topic consult, for example, Evans (2015). Let $\{f_{\theta}:\theta\in\Theta\}$ be a collection of densities on a
sample space $\mathcal{X}$ and $\pi$ be a prior on $\Theta.$ The posterior distribution of $\theta$ given that data $x$ is  $\pi(\theta\,|\,x)=\pi(\theta)f_{\theta}(x)/\int
_{\Theta}\pi(\theta)f_{\theta}(x)\,d\theta$. For an arbitrary parameter of interest $\psi=\Psi(\theta),$ the
prior and posterior densities of $\psi$ are denoted by $\pi_{\Psi}$ and $\pi_{\Psi}%
(\cdot\,|\,x),$ respectively. The relative belief ratio for a value $\psi$ is
then defined by $RB_{\Psi}(\psi\,|\,x)=\lim_{\delta\rightarrow0}\Pi_{\Psi
}(N_{\delta}(\psi\,)|\,x)/\Pi_{\Psi}(N_{\delta}(\psi\,))$, where $N_{\delta
}(\psi\,)$ is a sequence of neighbourhoods of $\psi$ converging   nicely (see, for example, Rudin (1974)) to
$\psi$ as $\delta\rightarrow0.$ Quit generally
\begin{equation}
RB_{\Psi}(\psi\,|\,x)=\pi_{\Psi}(\psi\,|\,x)/\pi_{\Psi}(\psi), \label{relbel}%
\end{equation}
the ratio of the posterior density to the prior density at $\psi.$ That is,
$RB_{\Psi}(\psi\,|\,x)$ is measuring how beliefs have changed that
$\psi$ is the true value from \textit{a priori} to \textit{a posteriori}. Note that, a relative
belief ratio is similar to a Bayes factor, as both are measures of evidence,
but the latter measures this via the change in an odds ratio.  A discussion about the relationship between relative belief ratios and Bayes factors is detailed in (Baskurt and Evans, 2013). In particular, when a Bayes factor is defined via a limit in the continuous
case, the limiting value is the corresponding relative belief ratio.

By a basic principle of evidence,  $RB_{\Psi}(\psi\,|\,x)>1$ means that the data led to an increase in the probability that $\psi$ is correct, and so there is
evidence in favour of $\psi,$ while $RB_{\Psi}(\psi\,|\,x)<1$ means that the data led
to a decrease in the probability that $\psi$ is correct, and so there is
evidence against $\psi,$. Clearly, when $RB_{\Psi}(\psi\,|\,x)=1$, then there is no
evidence either way.

Thus, the value $RB_{\Psi}(\psi_{0}\,|\,x)$  measures the evidence for the
hypothesis $\mathcal{H}_{0}=\{\theta:\Psi(\theta)=\psi_{0}\}.$ It is also important to calibrate whether this is strong
or weak evidence for or against $\mathcal{H}_{0}$. As suggested in Evans
(2015), a useful calibration of
$RB_{\Psi}(\psi_{0}\,|\,x)$ is obtained by computing the tail probability
\begin{equation}
\Pi_{\Psi}(RB_{\Psi}(\psi\,|\,x)\leq RB_{\Psi}(\psi_{0}\,|\,x)\,|\,x).
\label{strength}%
\end{equation}
One way to view (\ref{strength}) is as the posterior probability that the true value of $\psi$ has a relative
belief ratio no greater than that of the hypothesized value $\psi_{0}.$ When $RB_{\Psi}(\psi_{0}\,|\,x)<1,$ so there is evidence
against $\psi_{0},$ then a small value for (\ref{strength}) indicates a large
posterior probability that the true value has a relative belief ratio greater
than $RB_{\Psi}(\psi_{0}\,|\,x)$ and so there is strong evidence against
$\psi_{0}.$ When $RB_{\Psi}(\psi_{0}\,|\,x)>1,$ so there is evidence in favour
of $\psi_{0},$ then a large value for (\ref{strength}) indicates a small
posterior probability that the true value has a relative belief ratio greater
than $RB_{\Psi}(\psi_{0}\,|\,x))$ and so there is strong evidence in favour of
$\psi_{0},$ while a small value of (\ref{strength}) only indicates weak
evidence in favour of $\psi_{0}.$

\section{The Dirichlet Process}

In this section, a concise summary of the Dirichlet process is given. Because of its attractive features, the Dirichlet process, formally introduced in Ferguson (1973), is considered the most well-known and widely used prior in Bayesian nonparametric inference. Consider a space $\mathfrak{X}$ with a $\sigma-$algebra $\mathcal{A}$ of subsets of $\mathfrak{X}$. Let $H$ be a fixed probability measure on $(\mathfrak{X},\mathcal{A})$, called the \emph{base measure},  %
 and $a$ be a positive number, called the \emph{concentration parameter}. Following Ferguson (1973), a
random probability measure $P=\left\{  P(A)\right\}  _{A\in\mathcal{A}}$ is
called a Dirichlet process on $(\mathfrak{X},\mathcal{A})$ with parameters $a$
and $H$, denoted by $DP(a,H)$, if for any finite measurable partition $\{A_{1},\ldots,A_{k}\}$ of
$\mathfrak{X}$ with $k \ge 2$, $\left(  P(A_{1}%
),\ldots\,P(A_{k})\right)\sim \text{Dirichlet}(aH(A_{1}),\ldots,$ $aH(A_{k}))$. It is assumed that if
$H(A_{j})=0$, then $P(A_{j})=0$ with a probability one. Note that, for any
$A\in\mathcal{A},$ $P(A) \sim \text{Beta}(aH(A),(1-H(A))$ and so ${E}(P(A))=H(A)\ \ $ and ${Var}(P(A))=H(A)(1-H(A))/(1+a).$ Thus, $G$ can be viewed as the center of the process. On the other hand, $a$ controls concentration, as the larger value of $a$, the more likely that $P$ will be close to $G$. We refer the reader to Al-Labadi and Abdelrazeq (2017) for additional interesting asymptotic properties of the Dirichlet process and other nonparametric priors.

A distinctive feature of the Dirichlet process, among many other nonparametric priors,  is its conjugacy property. Specifically, if
$x=(x_{1},\ldots,x_{n})$ is a sample from $P\sim DP(a,H)$, then the posterior
distribution of $P$ is $P\,|\,x=P_x\sim DP(a+n,H_{x})$ where
\begin{equation}
H_{x}=a(a+n)^{-1}H+n(a+n)^{-1}F_{n}, \label{DP_posterior}%
\end{equation}
with $F_{n}=n^{-1}\sum_{i=1}^{n}\delta_{{x}_{i}}$ and $\delta_{x_{i}}$ is the
Dirac measure at $x_{i}.$ Notice that, $H_{x}$
is a convex combination of the prior base distribution and the empirical
distribution. Clearly,   $H_{x}\to H$  as $a \to \infty$ while   $H_{x}\to F_n$ as $a \to 0$.

Following Ferguson (1973), $P\sim{DP}(a,H)$ has the following series representation
\begin{equation}
P=\sum_{i=1}^{\infty}J_{i} \delta_{Y_{i}}, \label{series-dp}%
\end{equation}
where $\Gamma_{i}=E_{1}+\cdots+E_{i}$, $E_{i} \overset{i.i.d.}\sim \text{exponential}(1)$, $Y_{i} \overset{i.i.d.}\sim  H$
 independent of $\Gamma_{i}$, $L(x)=a\int_{x}^{\infty}t^{-1}e^{-t}dt,x>0,$ $L^{-1}(y)=\inf
\{x>0:L(x)\geq y\}$ and $J_{i}=L^{-1}(\Gamma_{i})/\sum_{i=1}^{\infty
}{L^{-1}(\Gamma_{i})}$.  It follows clearly from (\ref{series-dp}) that a
realization of the Dirichlet process is a discrete probability measure.  This
is true even when the base measure is absolutely continuous. One could resemble the discreteness of $P$ with the discreteness of $F_n$. Note that, since data is always measured to finite accuracy, the true distribution being sampled from is discrete. This makes the  discreteness property of $P$ with no practical significant limitation.  Indeed, by imposing the weak topology, the support for the Dirichlet process
is quite large. Specifically, the support for the Dirichlet process is the set
of all probability measures whose support is contained in the support of the
base measure. This means if the support of the base measure is $\mathfrak{X}$,
then the space of all probability measures is the support of the Dirichlet
process. In particular, if we have a normal base measure, then the Dirichlet
process can choose any probability measure.

Zarepour and Al-Labadi (2012) derived the following  series
approximation with monotonically decreasing weights for the Dirichlet process
\begin{equation}
P_{N}=\sum_{i=1}^{N}J_{i}\delta_{Y_{i}}, \label{eq11}
\end{equation}
where  $Y_{i}$  and $\Gamma
_{i}$ are as defined in (\ref{series-dp}), $G_{a/N}$ be the co-cdf of the g$\text{amma}(a/N,1)$
distribution and $J_{i}={G_{a/N}^{-1}(\Gamma_{i}/\Gamma_{N+1})/}\sum
_{j=1}^{N}{G_{a/N}^{-1}(\Gamma_{j}/\Gamma_{N+1})}$. They proved that, as $N \to
\infty$, $P_{N}$ converges almost surely to (\ref{series-dp}).  Note that ${G_{a/N}^{-1}(p)}$ is the $(1-p)$-th quantile of the
g$\text{amma}(a/N,1)$ distribution. This provides the following
algorithm. $\bigskip$

\noindent\textbf{Algorithm A: Approximately generating a value from }%
$DP(a,H)$\textbf{\smallskip}

\noindent1. Fix a relatively large positive integer $N$.\textbf{\smallskip}

\noindent2. For $i=1,\ldots,N$, generate  $Y_{i}\overset{i.i.d.} \sim H$. %
\textbf{\smallskip}

\noindent3. Independent of $\left(  Y_{i}\right)  _{1\leq
i\leq N}$,  for $i=1,\ldots,N+1,$ generate  $E_{i}\overset{i.i.d.} \sim\,$%
exponential$(1)$ and put $\Gamma_{i}=E_{1}+\cdots+E_{i}.$\textbf{\smallskip}

\noindent4. For $i=1,\ldots,N$, compute $G_{a/N}^{-1}\left(  {\Gamma_{i}%
}/{\Gamma_{N+1}}\right)  .$\textbf{\smallskip}

\noindent5. Use $P_N$ in (\ref{eq11}) to obtain an approximate value from
$DP(a,H)$.$\bigskip$

\noindent  For other simulation methods for the Dirichlet process, see, for instance,
Bondesson (1982), Sethuraman (1994), Wolpert and Ickstadt (1998) and Al-Labadi and Zarepour (2014b).

Throughout the paper, the notation $P$ could refer to either a
probability measure or its corresponding cdf where the context determines the
appropriate interpretation. That is, $P((-\infty,t])=P(t)$ for  all $t \in \mathbb{R}$.

\section{Cram\'{e}r-von Mises Distance}

A well-known and widely used distance between two distributions is the Cram\'{e}r-von Mises
Distance.  For cdf's $F$ and $G$ this is defined as $$d_{CvM}(F,G)=\int
_{-\infty}^{\infty}\left(  F(x)-G(x)\right)  ^{2}G(dx).$$

The next lemma demonstrates that, as sample sizes get large, the Cram\'{e}r-von Mises distance between posterior distributions of Dirichlet processes converges to the Cram\'{e}r-von Mises distance  between the true  distributions generated the data.

\begin{lemma}
\label{BSP3}
  Given two independent samples $x=(x_1,\ldots,x_{n_1}) \overset {i.i.d.} \sim F$ and $y=(y_1,\ldots,y_{n_2}) \overset {i.i.d.} \sim G$, with $F$ and $G$ being continuous  cdf's. Let $P\sim DP(a_1,H_1)$, $Q\sim DP(a_2,H_2)$, $P|x=P_x$ and $Q|y=Q_y$. Then, as $n_1,n_2\to \infty$, $d(P_{x},Q_{y}) \overset{a.s.}\to d(F,G)$.
\end{lemma}

The next corollary shows that the posterior distribution of
$d_{CvM}(P_x,Q_y)$ becomes concentrated around 0 as sample sizes
increase if and only if $\mathcal{H}_{0}$ holds. The proof follows straightforwardly from Lemma \ref{BSP3}.
\begin{lemma}
\label{cvm4} Let $x=(x_1,\ldots,x_{n_1}) \overset {i.i.d.} \sim F$ and $y=(y_1,\ldots,y_{n_2}) \overset {i.i.d.} \sim G$, with $F$ and $G$ being continuous  cdf's. Let $P\sim DP(a_1,H_1)$ and $Q\sim DP(a_2,H_2)$. As $n_1,n_2\to \infty$, (i) if
$\mathcal{H}_{0}$ is true, then $d_{CvM}\left(  P_{x},Q_{y}\right)
\overset{a.s.}{\rightarrow}0$ and (ii) if $\mathcal{H}_{0}$\ is false, then
$\lim\inf d_{CvM}(P_{x},Q_{y})\overset{a.s.}{>}0.$
\end{lemma}

The following result allows the use of the approximation (\ref{eq11}) when considering the prior and posterior distributions of the
Cram\'{e}r-von Mises distance.

\begin{lemma}
\label{cvm3} Let  $P\sim DP(a_1,H_1)$ and $Q\sim DP(a_2,H_2)$. Let  $P_{N_1}$  and $Q_{N_2}$ be two approximations of $P$ and $Q$, respectively, as defined in (\ref{eq11}). Then, as $N_1,N_2 \to \infty$,
$d_{CvM}\left(  P_{N_1},Q_{N_2}\right)  \overset{a.s.}{\rightarrow}d_{CvM}\left(
P,Q\right).$
\end{lemma}

The next lemma demonstrates that the distribution of the distance between two Dirichlet processes is independent from the base measures. This result will play a key role in the proposed approach.

\begin{lemma} \label{BSP4}
Let $P \sim DP(a_1,H_1)$ and $Q\sim DP(a_1,H_2)$, where $H_1$ and $H_2$ are continuous. If $H_1=H_2$, then the distribution of $d_{CvM}\left(P,Q\right)$ does not depend on  $H_1$ and $H_2$.
\end{lemma}


\section{The Approach}
Let  $x=(x_1,\ldots,x_{n_1}) \overset {i.i.d.} \sim F$ and $y=(y_1,\ldots,y_{n_2}) \overset {i.i.d.} \sim G$  be independent samples with $F$ and $G$ being unknown continuous cdf's. The goal to test the null hypothesis $\mathcal{H}_0:~F=G$. To this end, we use the priors $P\sim DP(a_1,H_1)$ and $Q \sim DP(a_2,H_2)$ so, by (\ref{DP_posterior}), $P|x \sim DP(a_1+n_1,H_x)$ and $Q|y \sim DP(a_1+n_1,H_y)$.  From  Lemma \ref{BSP3}, $D_{x,y}=d_{CvM}(P_x, Q_y)$ almost surely approximate $D=d_{CvM}(F,G)$. Thus, it looks clear that if $\mathcal{H}_{0}$ is true, then the
posterior distribution of the distance between $P$  and $Q$ should be more concentrated about $0$ than the prior distribution of the distance between $P$ and $Q.$ For example, in Figure 1-a (see Example 1), since $\mathcal{H}_{0}$ is true,
the plot of the posterior density of $D_{x,y}$ is much more concentrated about 0 than the
the plot of the prior density of $D$. So, the proposed test includes a comparison of the concentrations of the
prior and posterior distributions of $d_{CvM}$ via a relative belief ratio
based on $d_{CvM}$ with the interpretation as discussed in$\ $Section 2.

The  success of the approach depends significantly on a suitable selection of the parameters of $DP(a_1,H_1)$ and $DP(a_2,H_2)$. As illustrated below, inappropriate values of the parameters can lead to a failure in computing $d_{CvM}$. We discuss first setting  values of $H_1$ and $H_2$. By Lemma \ref{BSP4},  the distribution of  $d_{CvM}\left(P,Q\right)$ is  independent from the choice of the base measures when $H_1=H_2$, where both need to be continuous. Thus,  we suggest to set $H_1=H_2=N(0,1)$,  although other choices of continuous distributions are certainly possible. An additional and important reason supporting the choice of $H_1=H_2$ is to avoid prior-data conflict  (Evans and Moshonov, 2006; Al-Labadi and Evans, 2017).  Prior-data conflict means that there is a tiny overlap between the effective support regions of $DP(a_1,H_1)$  and $DP(a_2,H_2)$. In this context, the existence of prior-data conflict can yield to  a failure in computing the distribution of $d_{CvM}\left(  P,Q\right)  $ about 0. To avoid prior-data conflict, it is necessary that $H_1$ and $H_2$ share the same effective support (note that, $P$ and $Q$ have the same  support as $H_1$ and $H_2$, respectively), which can  certainly be secured by setting  $H_1=H_2$. The effect of prior-data conflict is demonstrated in Section 7, Table \ref{tab2}.

The selection of $a_1$ and $a_2$ is also important.  It is possible to consider several values of $a_1$ and $a_2$. In general, the values of $a_1$ and $a_2$ depends in $n_1$ and $n_2$, respectively. As indicated in Al-Labadi and Zarepour (2017), $a_i$ should be chosen to have a value at most $0.5n_i,  i=1,2$ as otherwise the prior may become too influential. Holmes  et al. (2015) recommend using values between 1 and 10 and checking the sensitivity of the results to the chosen values. The following algorithm outlines a procedure for selecting the concentration parameters. \smallskip

\noindent\textbf{Algorithm B: Selection of concentration parameters\smallskip}

\noindent1. Start by setting $a_1=a_2=1$ and compute the relative belief ratio and its strength. Algorithm C in the next section addresses such computations.\textbf{\smallskip}

\noindent 2. Consider more concentrated priors by setting larger values of $a_1$ and $a_2$.

\noindent 3. Compute the corresponding  relative belief ratio. There are two scenarios:

\begin{enumerate}
\item [a.] If the value of the relative belief ratio in step 1 is less (greater) than 1 and  the new value is less (greater) than 1, then there is an evidence against (in favour) $\mathcal{H}_{0}$.\textbf{\smallskip}

\item [b.] If  the value of the relative belief ratio in step 1 is greater than 1 and the new value is  greater (less) than 1, then this is an evidence against (in favour) $\mathcal{H}_{0}$.\textbf{\smallskip}
\end{enumerate}

Algorithm B is further explored in  Table \ref{tab1} of Section 7. In most cases, setting $a_1=a_2=1$ is found to be adequate.  Holmes  et al. (2015) recommend using values between 1 and 10 and checking the sensitivity of the results to the chosen values.

\section{Computations}

Closed forms of the densities of $D=d_{CvM}(P,Q)$ and  $D_{x,y}=d_{CvM}(P_x,Q_y)$ are typically not available. Thus, the relative belief
ratios need to be approximated via simulation. The following gives a computational algorithm to test $\mathcal{H}_{0}$. This algorithm is a revised version of Algorithm B of Al Labadi and Evans (2018).
\noindent\textbf{Algorithm C: Relative belief algorithm for the two-sample problem\smallskip}

\noindent1. Use Algorithm A to (approximately) generate  a $P$  from
$DP(a_1=1,N(0,1))$ and a $Q$ from
$DP(a_2=1,N(0,1))$. \textbf{\smallskip}

\noindent2. Compute $d_{CvM}(P,Q)$.\textbf{\smallskip}

\noindent3. Repeat steps (1)-(2) to obtain a sample of $r_{1}$ values from the
prior of $D$.\textbf{\smallskip}

\noindent4. Use Algorithm A to (approximately) generate a $P_x$ from
$DP(1+n_1,H_{x})$ and  $Q_x$ from
$DP(1+n_2,H_{y})$.\textbf{\smallskip}

\noindent5. Compute $d_{CvM}(P_x,Q_y)$.\textbf{\smallskip}

\noindent6. Repeat steps (4)-(5) to obtain a sample of $r_{2}$ values  of $D_{x,y}$.\textbf{\smallskip}

\noindent7. Let $M$ be a positive number. Let $\hat{F}_{D}$ denote the
empirical cdf of $D$ based on the prior sample in (3) and for $i=0,\ldots,M,$
let $\hat{d}_{i/M}$ be the estimate of $d_{i/M},$ the $(i/M)$-th prior
quantile of $D.$ Here $\hat{d}_{0}=0$, and $\hat{d}_{1}$ is the largest value
of $d$. Let $\hat{F}_{D}(\cdot\,|\,x,y)$ denote the empirical cdf of $D$ based
on the posterior sample in 6. For $d\in\lbrack\hat{d}_{i/M},\hat
{d}_{(i+1)/M})$, estimate $RB_{D}(d\,|\,x,y)={\pi_D(d|x,y)}/{\pi_D(d)}$ by
\begin{equation}
\widehat{RB}_{D}(d\,|\,x,y)=M\{\hat{F}_{D}(\hat{d}_{(i+1)/M}\,|\,x,y)-\hat{F}%
_{D}(\hat{d}_{i/M}\,|\,x,y)\}, \label{rbest}%
\end{equation}
the ratio of the estimates of the posterior and prior contents of $[\hat
{d}_{i/M},\hat{d}_{(i+1)/M}).$ It follows that,   we estimate $RB_{D}(0\,|\,x,y)={\pi_D(0|x,y)}/{\pi_D(0)}$
 by
$\widehat{RB}_{D}(0\,|\,x,y)=$ $M\widehat{F}_{D}(\hat{d}_{p_{0}}\,|\,x,y)$ where
$p_{0}=i_{0}/M$ and $i_{0}$ is chosen so that $i_{0}/M$ is not too small
(typically $i_{0}/M\approx0.05)$.\textbf{\smallskip}

\noindent8. Estimate the strength $DP_{D}(RB_{D}(d\,|\,x,y)\leq RB_{D}%
(0\,|\,x,y)\,|\,x,y)$ by the finite sum
\begin{equation}
\sum_{\{i\geq i_{0}:\widehat{RB}_{D}(\hat{d}_{i/M}\,|\,x,y)\leq\widehat{RB}%
_{D}(0\,|\,x,y)\}}(\hat{F}_{D}(\hat{d}_{(i+1)/M}\,|\,x,y)-\hat{F}_{D}(\hat
{d}_{i/M}\,|\,x,y)). \label{strest}%
\end{equation}

\noindent For fixed $M,$ as $r_{1}\rightarrow\infty,r_{2}\rightarrow\infty,$
then $\hat{d}_{i/M}$ converges almost surely to $d_{i/M}$ and (\ref{rbest})
and (\ref{strest}) converge almost surely to $RB_{D}(d\,|\,x)$ and
$DP_{D}(RB_{D}(d\,|\,x,y)\leq RB_{D}(0\,|\,x,y)\,|\,x,y)$, respectively.

\noindent9. As detailed in  Algorithm B,  repeat steps (1)-(8) for larger values of $a_1$ and $a_2$.

The following proposition establishes the consistency of the approach to the two-sample problem
as sample size increases. So the procedure performs correctly as sample size increases when
$\mathcal{H}_{0}$ is true. The proof follows immediately from Evans (2015), Section 4.7.1.
\begin{proposition}
\label{cvm6}Consider the discretization $\{[0,d_{i_{0}/M}),[d_{i_{0}%
/M},d_{(i_{0}+1)/M}),\ldots,$\newline$[d_{(M-1)/M},\infty)\}$. As
$n_1,n_2\rightarrow\infty,$ (i) if $\mathcal{H}_{0}$ is true, then
\begin{align*}
&  RB_{D}([0,d_{i_{0}/M})\,|\,x,y)\overset{a.s.}{\rightarrow}1/DP_{D}%
([0,d_{i_{0}/M})),\\
&  RB_{D}([d_{i/M},d_{(i+1)/M})\,|\,x,y)\overset{a.s.}{\rightarrow}0\text{
whenever }i\geq i_{0},\\
&  DP_{D}(RB_{D}(d\,|\,x,y)\leq RB_{D}(0\,|\,x,y)\,|\,x,y)\overset{a.s.}%
{\rightarrow}1,
\end{align*}
and (ii) if $\mathcal{H}_{0}$ is false and $d_{CvM}(P, Q)\geq d_{i_{0}/M}$, then $RB_{D}([0,d_{i_{0}/M})\,|\,x,y)\overset
{a.s.}{\rightarrow}0$ and $DP_{D}(RB_{D}(d\,|\,x,y)\leq RB_{D}%
(0\,|\,x,y)\,|\,x,y)\overset{a.s.}{\rightarrow}0.$
\end{proposition}

\section{Examples}
In this section, the approach is illustrated through three examples. In Examples 1 and 2,  the methodology is assessed using simulated samples from
a variety of distributions and in Example 3 an application to a real data set
is presented.

The following notation is used for the distributions in the
tables, namely, $N(\mu,\sigma)$ is the normal distribution with mean $\mu$ and standard deviation $\sigma$, $t_r$ is the $t$ distribution with $r$ degrees of freedom, exp$(\lambda)$ is the exponential distribution with mean $\lambda$ and  $U(a,b)$ is the uniform distribution over $[a,b]$.  For all cases, we set  $N_1=N_2=1000$ in Algorithm A and $r_1=r_2=2000$, $M=20$ in Algorithm B. The results are also compared with the frequentist Cram\'er-von Mises (CvM) test. To calculate p-values of the CvM test, the \textbf{\textsf{R}} function ``cramer.test" is used. We also compared our results with the Bayesian nonparametric tests of Holmes  et al. (2015) and Al-Labadi and Zarepour (2017). Since the obtained results are similar in these tests, we reported only the results of the new approach. \smallskip

\noindent {\textbf{Example 1.}}  \label{example1}
Consider samples generated from the  distributions in Table \ref{tab1}, where each sample is of size 50 (Case 1- Case 9). These distributions are also considered in Holmes et al. (2015) and Al-Labadi and Zarepour (2017). To study the sensitivity of the approach  to the choice of concentration parameters,   various values of $a_1$  and $a_1$  are considered. The results are reported in Table \ref{tab1}. Recall that, we want $RB>1$ and the strength close to 1 when
$\mathcal{H}_{0}$ is true and $RB<1$ and the strength close to 0 when
$\mathcal{H}_{0}$ is false. It follows that,  the methodology performs
perfectly in all cases. For example, in Case 1, since $RB=9.4$ and strength$=1$, there is no reason to doubt that the two sampling distributions are not identical. On the other hand, in Case 2, since $RB=0$ and strength$=0$, the two samples are drawn from two different distributions. We point out that the standard Cram\'er-von Mises test failed to recognize the difference in Case 6 (i.e., $x\sim N(0,1)$ and $y\sim t_{0.5}$). Notice that, in all cases, the appropriate
conclusion is attained with $a_1=a_2=1$. The other values of $a_1$ and $a_2$ considered in Table \ref{tab1} support the reached conclusions.

Figure 1 provides plots of the density of the prior distance and the posterior
distance for some cases in Example 1. It follows, for instance, from Figure 1 that the
posterior density of the distance is more concentrated about 0 than the prior
density of the distance when the two distributions are equal but not to the same degree otherwise.%

\begin{table}[htbp]
  \centering

    \begin{tabular}[c]{llccc}
    \hline

Samples&$a_1=a_2$ &$RB$(Strength) & p-values \\
    \hline
$x\sim N(0,1)$, $y\sim N(0,1)$&1&         9.40(1)&                               0.2977   \\
& 10&    8.54(1)&                                                \\
& 20&    4.48(0.776)&                                       \\\hline

$x\sim N(0,1)$, $y\sim N(1,1)$&1& 0(0)& 0.0000\\
& 10&    0(0)&                                      \\
& 20&    0(0)&                                      \\\hline

$x\sim N(0,1)$, $y\sim N(0,2)$&1& 0(0)& 0.0030  \\
& 10&    0.08(0.004)&                                       \\
& 20&    0(0)&                                      \\\hline

$x\sim N(0,1)$, $y\sim 0.5N(-2,1)+0.5N(2,1)$&1& 0(0)& 0.0000 \\
& 10&    0(0)&                                      \\
& 20&    0(0)&                                      \\\hline

$x\sim N(0,1)$, $y\sim t_{3}$&1&9.40(1)& 0.4316\\
& 10&    8.60(1)&                                                \\
& 20&    5.78(1)&                                                \\\hline

$x\sim N(0,1)$, $y\sim t_{0.5}$&1&0(0) & 0.1169 \\
& 10&    0.28(0.023)&                                       \\
& 20&    0.04(0.002)&                                       \\\hline

$\log x\sim N(0,1)$,  $\log y\sim N(1,1)$&1&0.02(0.001)& 0.0000               \\
& 10&    0.02(0.001)&                                       \\
& 20&    0.02(0.001)&                                       \\\hline

$x\sim \text{exp}(1)$,  $ y\sim \text{exp}(2)$&1&0.10(0) &  0.0020         \\
& 10&    0.12(0.006)&                                       \\
& 20&    0.06(0.003)&                                       \\\hline

$x\sim \text{exp}(1)$, $ y\sim \text{exp}(1)$&1&9.02(1) & 0.6134           \\
& 10&    7.30(1)&                                                \\
& 20&    4.86(1)&                                                \\
\hline
     \end{tabular}%
 \caption{Relative belief ratios and strengths for testing the equality of the two distributions generating the samples in Example 1. p-values of the (frequentist) Cram\'er-von Mises test are also reported.} \label{tab1}%
\end{table}%

\begin{figure}[htbp]
\centering
\subfigure [$x\sim N(0,1)$, $y\sim N(0,1)$] {\epsfig{figure=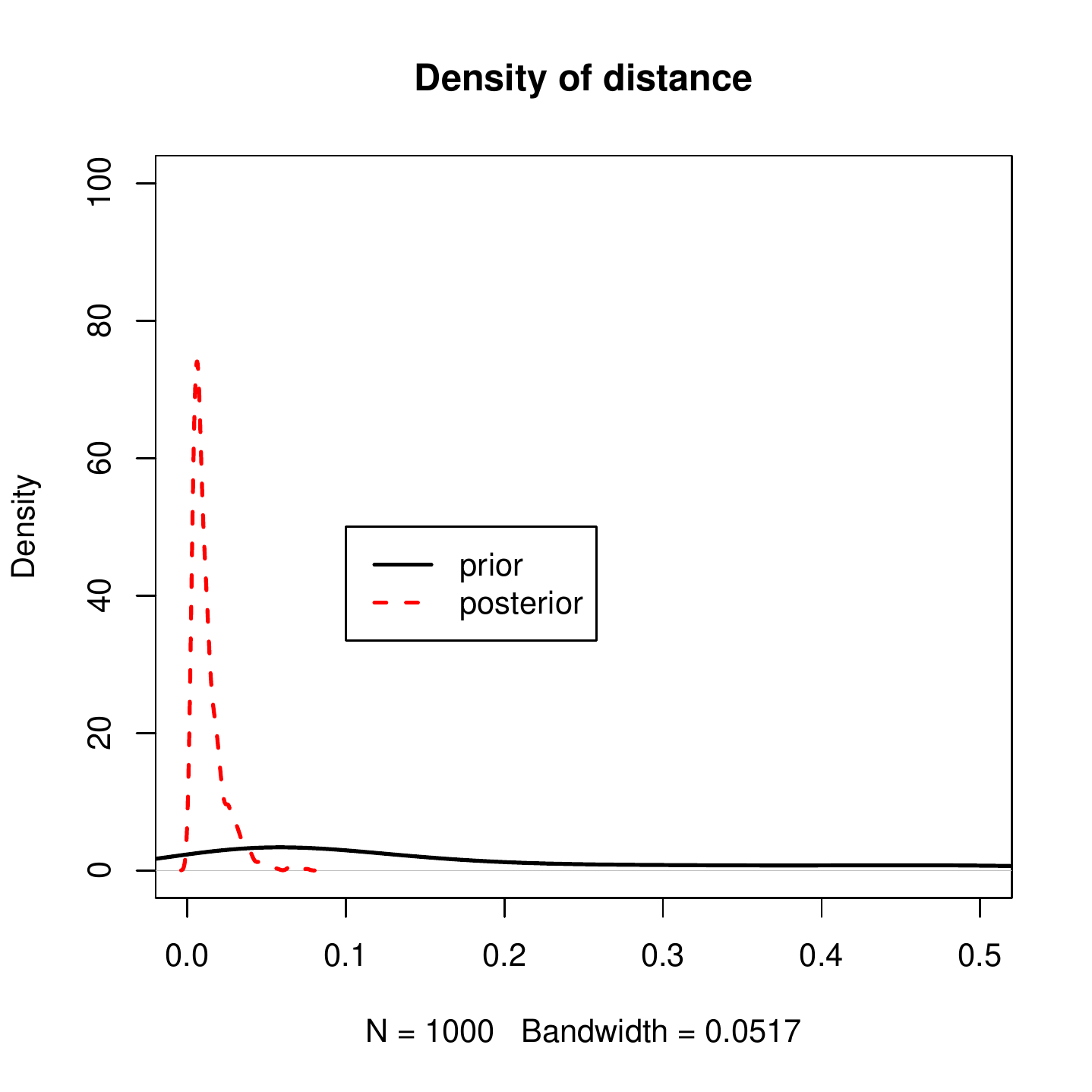,width=2.3in}}
\subfigure[$x\sim N(0,1)$, $y\sim N(1,1)$]{\epsfig{figure=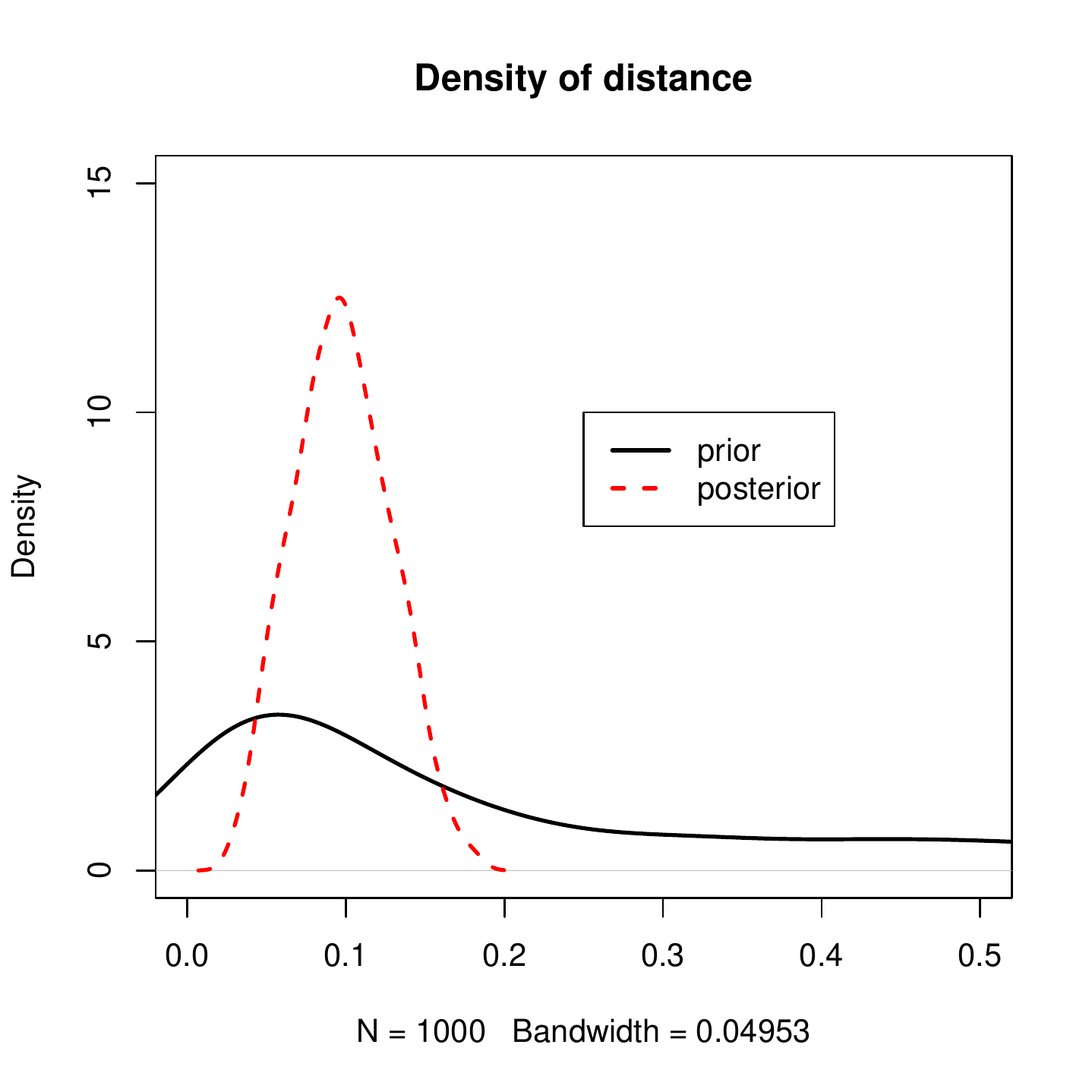,width=2.3in}}
\par
\subfigure [$x\sim \text{exp}(1)$, $y\sim \text{exp}(2)$] {\epsfig{figure=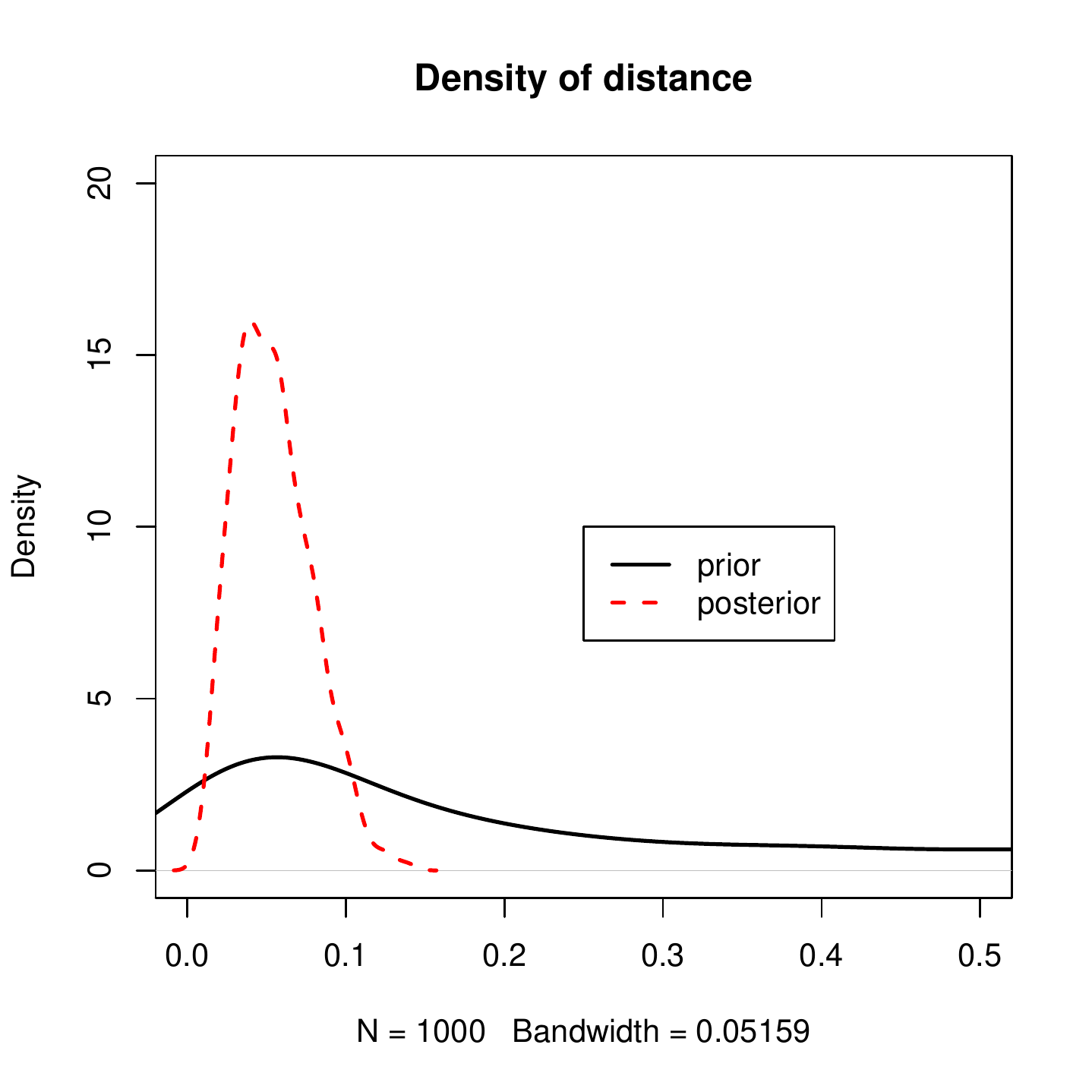,width=2.3in}}
\subfigure[$x\sim \text{exp}(1)$, $y\sim \text{exp}(1)$]{\epsfig{figure=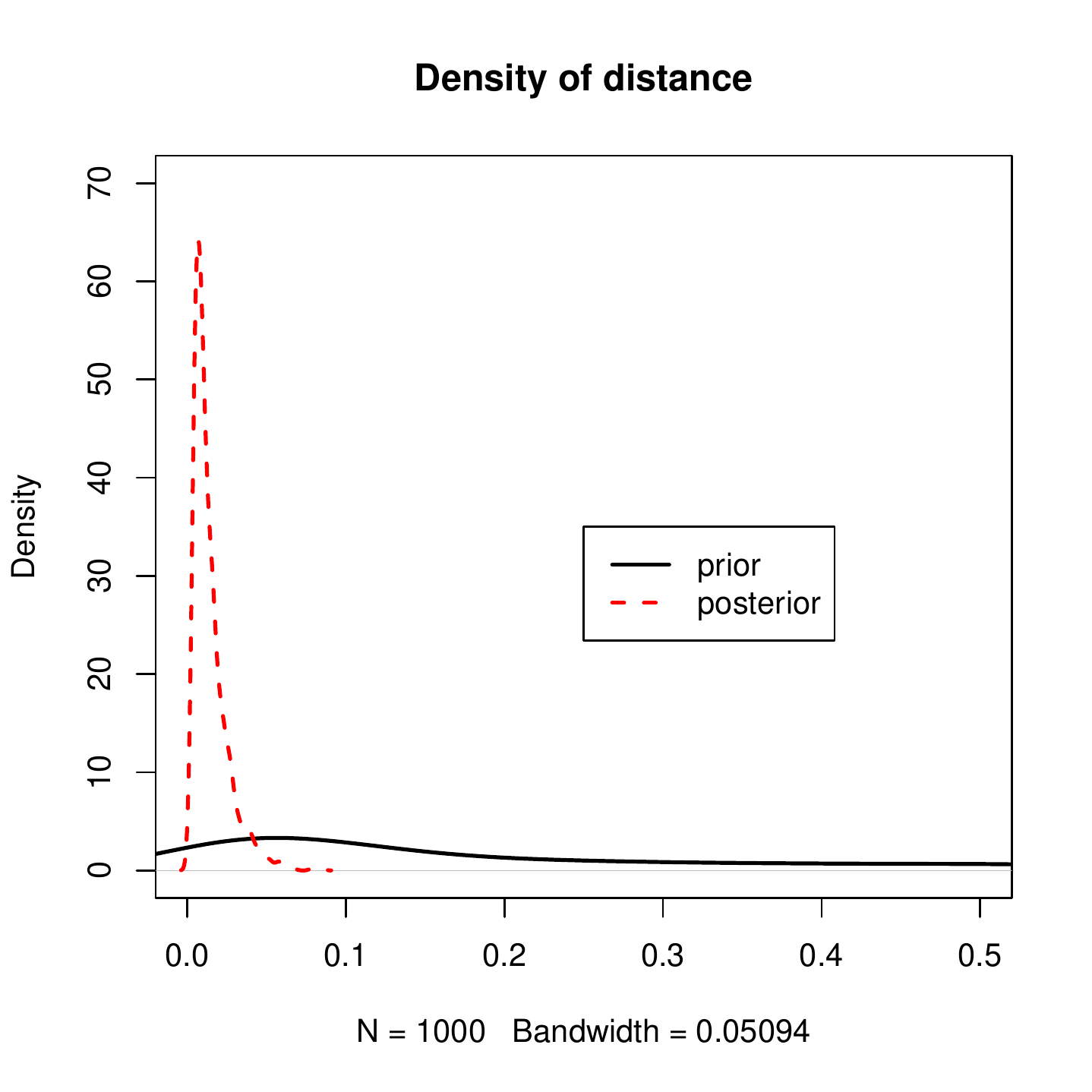,width=2.3in}}\caption{Plots of prior density versus posterior density of distance for some cases in Table 1.}%
\label{fig:SubF1}%
\end{figure}

It is also interesting to consider the effect of prior-data conflict on the
methodology. As discussed in Section 5, prior-data
conflict will occur whenever there is only a tiny overlap between  $H_1$ and $H_2$. Table \ref{tab2} gives the outcomes when  $x\sim N(0,1)$ and $y\sim N(1,1)$ for  a particular sample of sizes $n_1=n_2=50$ with various choices of $H_1$ and $H_2$.  Obviously,  only when $H_1=H_2$ we get the correct conclusion. This illustrates the importance of setting $H_1=H_2$ in the priors $DP(a_1,H_1)$  and $DP(a_1,H_1)$. \smallskip%

\begin{table}[htbp] \centering
\begin{tabular}
[c]{cllcc}\hline
Distribution & $H_1$ &$H_2$ &$RB$ (Strength)&p-value\\\hline
\multicolumn{1}{c}{$x\sim N(0,1)$, $y\sim N(1,1)$} & $N(0,1)$ & $N(0,1)$ &      $0\,(0)$&0.0000\\
\multicolumn{1}{c}{} & $N(-5,1)$ & $N(5,1)$ &      $20\,(1)$& \\
\multicolumn{1}{c}{} &$U(10,20)$  & $N(0,1)$  &      $20\,(1)$& \\
\multicolumn{1}{c}{} & $U(10,20)$ & $U(10,20)$ &      $0\,(0)$& \\
\hline
\end{tabular}
\caption{Study of prior-data conflict, where relative belief ratios and strengths with various choices of base measures $H_1$ and $H_2$ are computed.}\label{tab2}%
\end{table}%

Figure  2 also provides plots of the density of the prior distance and the posterior distance for the cases in Table \ref{tab2}. It follows that the correct conclusion is only obtained when $H_1=H_2$.

\begin{figure}[htbp]
\centering
\subfigure [$H_1=H_2=N(0,1)$] {\epsfig{figure=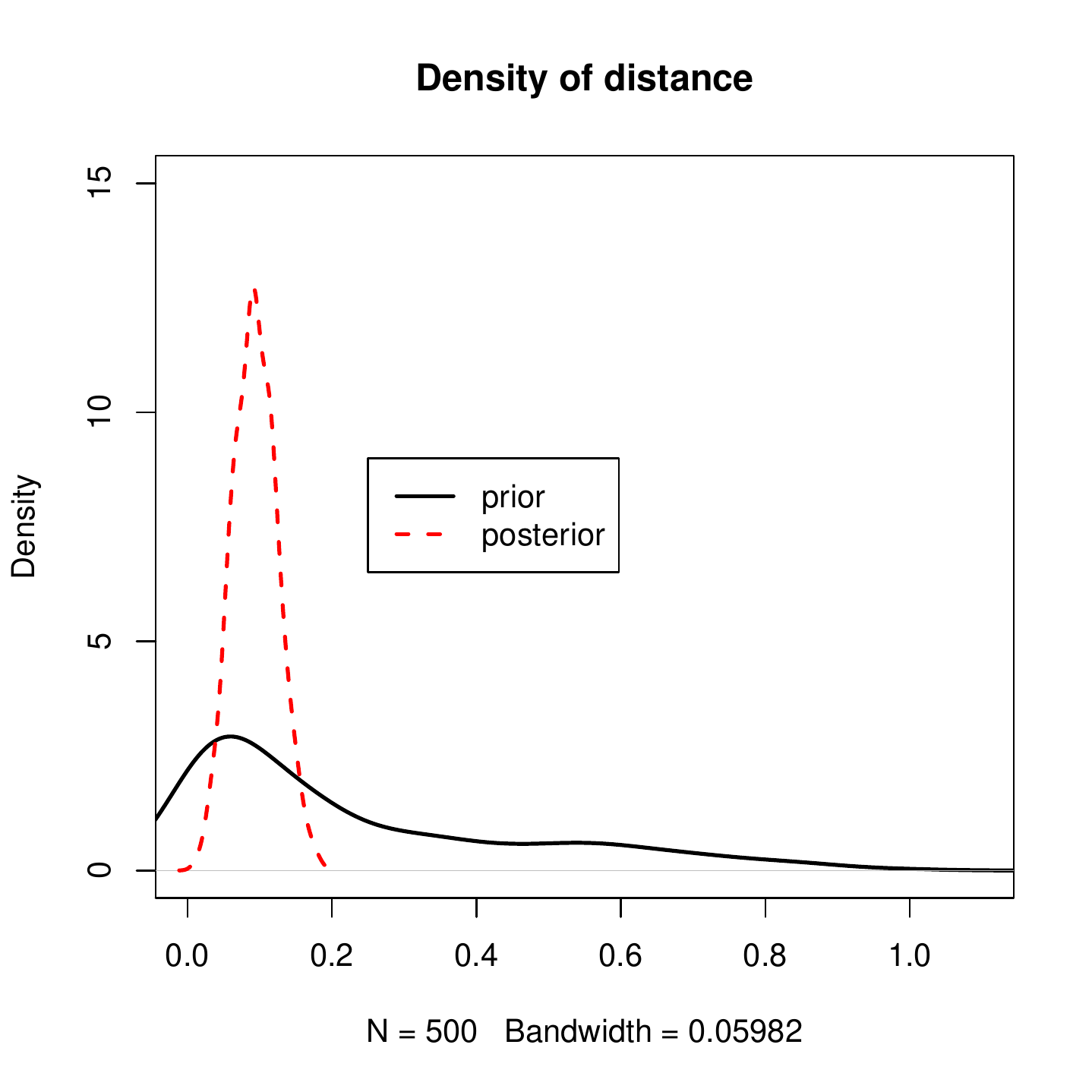,width=2.3in}}
\subfigure[$H_1=N(-5,1)$, $H_2=N(5,1)$]{\epsfig{figure=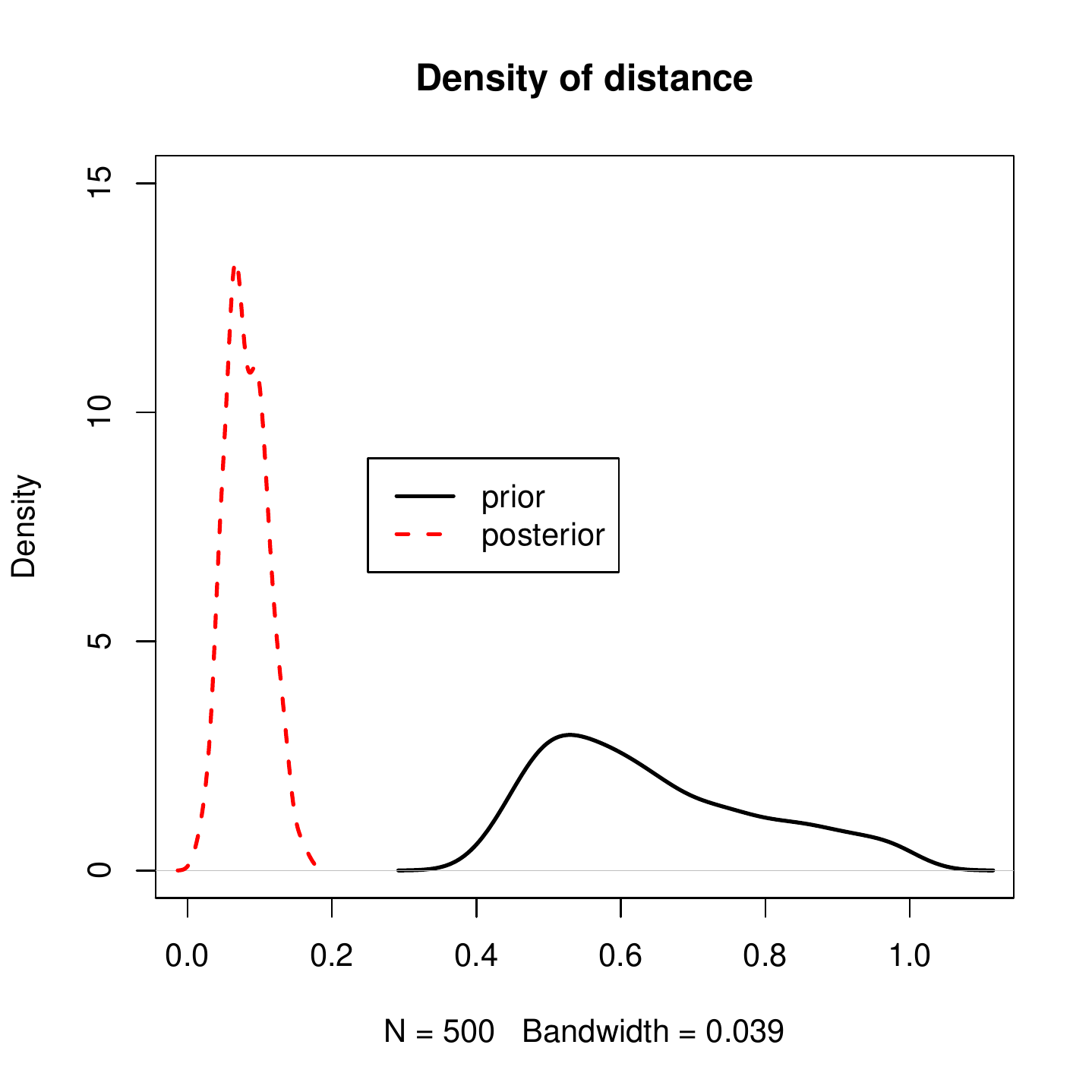,width=2.3in}}
\par
\subfigure [$H_1=U(10,20)$, $H_2=N(0,1)$] {\epsfig{figure=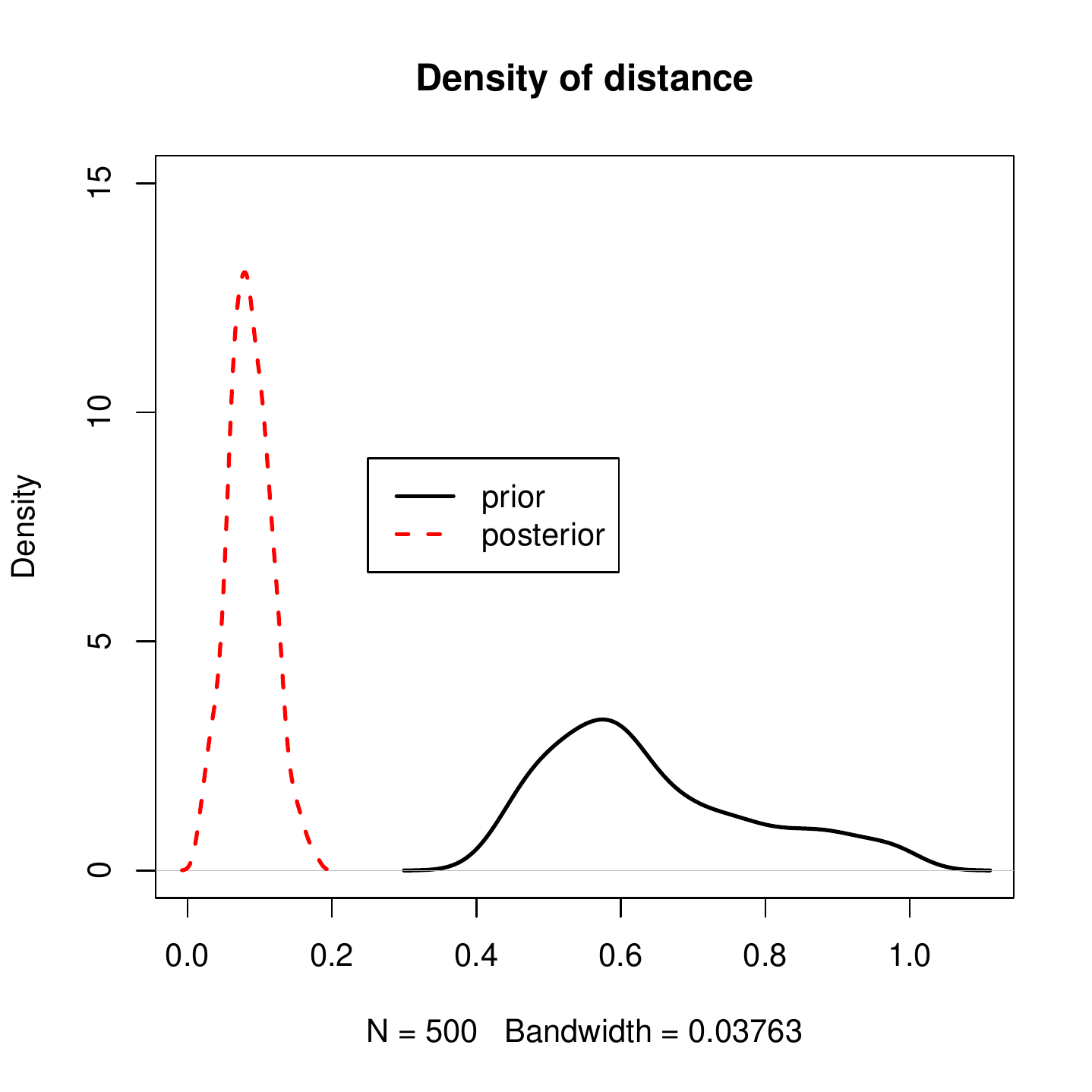,width=2.3in}}
\subfigure[$H_1=U(10,20)$, $H_2=U(10,20)$]{\epsfig{figure=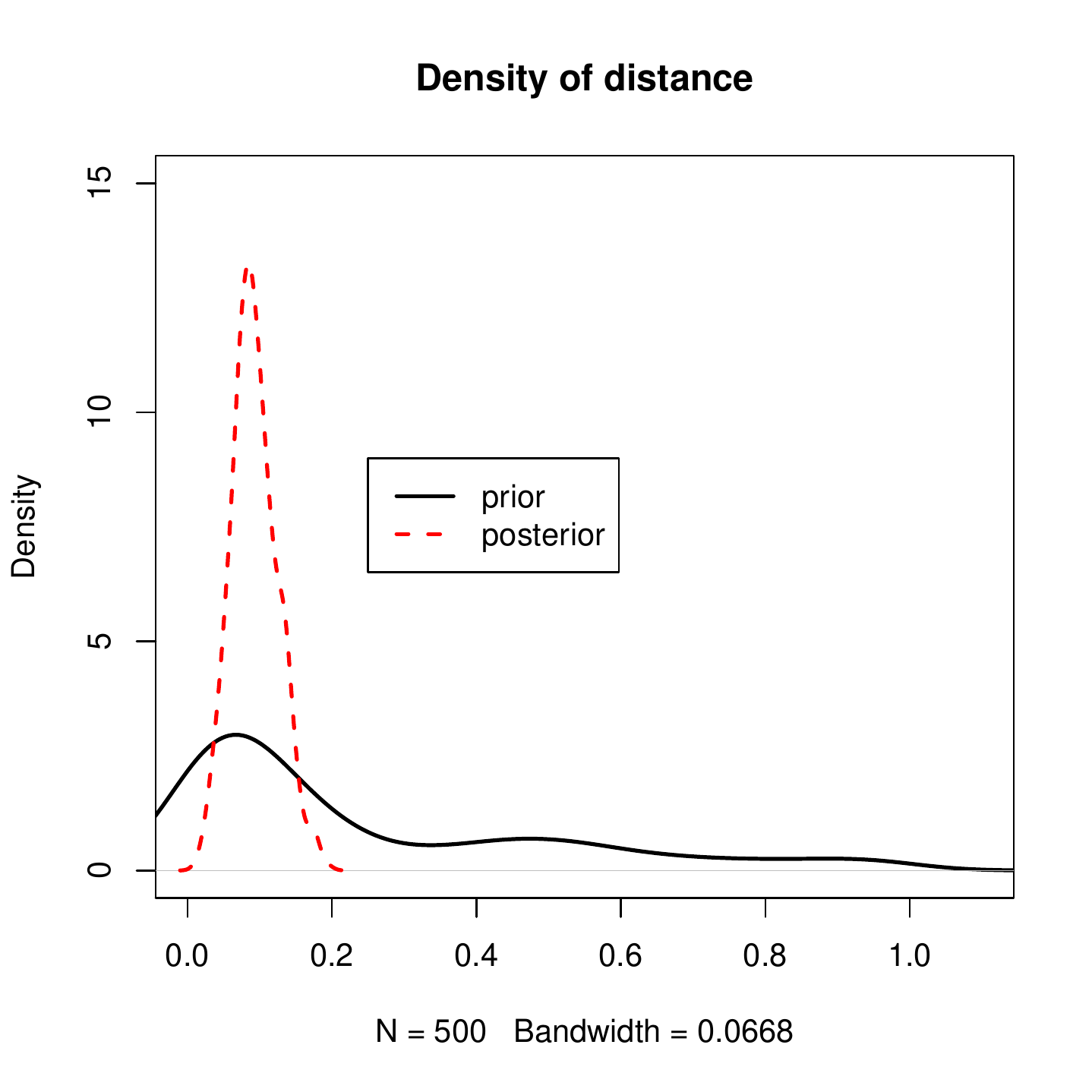,width=2.3in}}
\caption{Study of prior-data conflict, where plots of prior density versus posterior density of distance for the cases in Table \ref{tab2} are considered.}%
\label{fig:SubF1}%
\end{figure}

\noindent {\textbf{Example 2.}}  \label{example2} In this example, we explore the performance of the proposed test as  sample sizes increase. We consider samples from the distributions $x\sim N(0,1)$, $y\sim N(0,1)$ (Case 1) and  $x\sim N(0,1)$, $y\sim N(1,1)$ (Case 2). The results are summarized in Table \ref{tab3}. It follows  that the null hypothesis is not rejected in Case 1 but rejected in Case 2 . Clearly, the proposed approach  works well even with small sample sizes.

\begin{table}[htbp]
  \centering
      \begin{tabular}[c]{lccccc}
    \hline
 \multicolumn{1}{c}{\multirow{2}[4]{*}{Sample Sizes}} &\multicolumn{2}{c}{$x\sim N(0,1)$, $y\sim N(0,1)$} & \multicolumn{1}{c}{} &\multicolumn{2}{c}{$x\sim N(0,1)$, $y\sim N(1,1)$} \\  \cline{2-3} \cline{5-6}   &$RB$ (Strength)  & p-value   & & $RB$ (Strength) & p-value \\
    \hline
$n_1=n_2=5$                    &1.80(0.586)    &0.7083   &&0.36(0.02) &0.1628   \\
$n_1=n_2=10$  &1.24(0.250)& 0.8132   &&0.48(0.064) &0.1359    \\
$n_1=n_2=15$  &3.48(0.538)   & 0.9261   &&0.08(0.004)& 0.0069                 \\
$n_1=n_2=20$  &2.64(0.422)    &0.7103   && 0.12(0.010) &0.0170  \\
$n_1=n_2=30$  &5.60(1) &0.5864   &&0.08(0.006) & 0.0020           \\
$n_1=n_2=50$  &9.40(1)& 0.2977 &&0(0)              & 0.0030               \\
$n_1=n_2=100$                &13.08(1) &0.4236   &&0(0)         & 0.0000                           \\
$n_1=n_2=200$                &17.88(1)& 0.2697   &&0(0)         & 0.0000               \\ \hline
    \end{tabular}%
\caption{Relative belief ratios and strengths versus p-values.} \label{tab3}
\end{table}%

\noindent {\textbf{Example 3.}}  \label{example2} The proposed approach of the two-sample problems is illustrated on the chickwts data in \textbf{\textsf{R}}, where weights in grams are recorded for six groups of newly hatched chicks fed different supplements. The goal of this experiment was  to measure and compare the effectiveness of various feed supplements on the growth rate of chickens.  The first hypotheses of interest is to test whether the distributions  of weight of chicks  fed by soybean and linseed supplements differ. In the second hypothesis, we examine whether the distributions of  weight of chick for  sunflower and linseed groups differ. The ordered chick weights for the three samples are:

soybean: 158 171 193 199 230 243 248 248 250 267 271 316 327 329

linseed:  141 148 169 181 203 213 229 244 257 260 271 309

sunflower: 226 295 297 318 320 322 334 339 340 341 392 423

The values recorded in Table \ref{tab4}  do not support the evidence that the distributions  of the weight of chicks  fed by soybean and linseed supplements differ. On the other hand, they underline that the sunflower and linseed groups differ.

\begin{table}[htbp] \centering
\begin{tabular}
[c]{cccc}\hline
Samples&$a_1=a_2$ & $RB$ (Strength)& p-value\\\hline
$x$: soybean \& $y$: linseed &1 & $0.48\,(0.014)$&0.3487\\
&2 &  $2.50\,(0.717)$&\\
&3 &  $3.12\,(0.844)$&\\
&4 &  $ 3.14\,(0.843)$&\\
&5&  $3.34\,(0.833)$& \\\hline

$x$: soybean \& $y$: sunflower &1 & $0\,(0)$&0\\
&2 & $0\,(0)$&\\
&3 &  $0\,(0)$&\\
&4 &  $0\,(0)$&\\
&5 &  $0\,(0)$& \\\hline
\end{tabular}
\caption{Relative belief ratios and strengths for testing equality of  the distributions of chick weights for the soybean and linseed groups and the sunflower and linseed groups  of the chickwts data using various choices of
$a_1$ and $a_2$ in Example 3.}\label{tab4}%
\end{table}%

\section{Concluding Remarks}

A Bayesian approach for the two-sample problem  based on the use of the Dirichlet process and relative belief has been developed.  Implementing the approach is fairly simple and does not require obtaining a closed form of the relative belief ratio. Through several examples, it has been shown that the approach performs extremely well. While Cram\'{e}r-von Mises distance has been used in this paper, other distance
measures such as Anderson-Darling distance and the Kullback-Leibler distance are possible. It is also possible to extend the approach to cover the case of censored data.

\appendix

\section{Proofs}

\noindent \textbf{Proof of Lemma \ref{BSP3}}  For any cdf's
$P_{x}$ and $Q_{y}$, we have $d_{CvM}\left(  P_{x},Q_{y}\right)
=\int_{-\infty}^{\infty}\left(  P_{x}(z)-Q_{y}(z)\right)  ^{2}Q_{y}(dz)$. Since   $(P_{x}(z)-Q_{y}(z))^{2}\leq1$,
$P_{x}(z)\overset{a.s.}{\rightarrow}F(z)$  and $Q_{y}(z)\overset{a.s.}{\rightarrow}G(z)$ (James, 2008;  Al-Labadi and Abdelrazeq, 2017), the dominated convergence theorem completes the proof.
\endproof

\bigskip

\noindent \textbf{Proof of Lemma \ref{cvm3}} The proof is similar to the proof of Lemma 1. We include the proof for the sake of completeness.  For cdf's $P_{N_1}$ and $Q_{N_2}$, we have $d_{CvM}\left(P_{N_1},Q_{N_2}\right)
=\int_{-\infty}^{\infty}\left(  P_{N_1}(z)-Q_{N_2}(z)\right)  ^{2}$ $Q_{N_2}(dz)$. Since $(P_{N_1}(z)-Q_{N_2}(z))^{2}\leq1$,
$P_{N_1}(z)\overset{a.s.}{\rightarrow}P(z)$ and $Q_{N_2}(z)\overset{a.s.}{\rightarrow}Q(z)$, the result is followed by the dominating convergence theorem.
\endproof

\bigskip

\noindent \textbf{Proof of Lemma \ref{BSP4}}    
Since $H_1$ is nondecreasing, we have
 \begin{eqnarray*}
 \theta_i < t \ \ \text{if and only if} \ \ H_1(\theta_i) < H_1(t).
   \end{eqnarray*}
It follows from (\ref{series-dp}) that
 \begin{eqnarray*}
 P(t)=P\left((-\infty,t]\right)&=&\sum_{i=1}^{\infty} J_i \delta_{\theta_i}\left((-\infty,t]\right)=\sum_{i=1}^{\infty} J_i \delta_{H_1(\theta_i)}\left((0,H_1(t)]\right).
  \end{eqnarray*}
Observe that, since $(\theta_i)_{i \ge 1}$ is a sequence of
i.i.d. random variables with continuous distribution $H_1,$ for $i \ge 1,$ we have $U_i\overset{d}=H_1(\theta_i)$, where $\left(U_i\right)_{i \ge 1}$ is a sequence of i.i.d. random variables with a uniform distribution on $[0,1]$. Hence, $P(t)=P_{\lambda}(H_1(t)),$
where $P_{\lambda_1}\sim DP(a_1,\lambda)$ and $\lambda$ is the Lebesgue measure  on $[0,1]$. Similarly, $Q(t)=Q_{\lambda}(H_2(t)),$ where $Q_{\lambda}\sim DP(a_2,\lambda)$. Thus,
 \begin{eqnarray*}
d_{CvM}(P,Q)&=&\int
_{-\infty}^{\infty}\left(P(t)-Q(t)\right)  ^{2}Q(dt)\\
 &=&\int_{-\infty}^{\infty}\left(P_{\lambda}(H_1(t))-Q_{\lambda}(H_2(t))\right)^{2}Q_{\lambda}(H_2(dt))
  \end{eqnarray*}
If $H_1=H_2=H$, and since $H$ is continuous, we have
\begin{eqnarray*}
d_{CvM}(P,Q) &=&\int_{-\infty}^{\infty}\left(P_{\lambda}(H(t))-Q_{\lambda}(H(t))\right)^{2}Q_{\lambda}(H(dt))\\
  &=&\int_{0}^{1}\left(P_{\lambda}(z)-Q_{\lambda}(z)\right)^{2}Q_{\lambda}(dz).
  \end{eqnarray*}
  This shows that the distribution of $d_{CvM}(P,Q)$  does not depend on the base measures $H_1$ and $H_2$ whenever $H_1=H_2$.
\endproof

\end{document}